\documentclass[conference]{IEEEtran}
\IEEEoverridecommandlockouts
\usepackage{cite}
\usepackage{amsmath,amssymb,amsfonts}
\usepackage{algorithmic}
\usepackage{graphicx}
\usepackage{textcomp}
\usepackage{caption}
\usepackage{subcaption}
\usepackage{balance}
\usepackage{enumitem}
\usepackage{comment}

\usepackage[ruled,linesnumbered]{algorithm2e}
\usepackage{xcolor}
\def\BibTeX{{\rm B\kern-.05em{\sc i\kern-.025em b}\kern-.08em
    T\kern-.1667em\lower.7ex\hbox{E}\kern-.125emX}}
    
\newcommand{\R}{\mathbb{R}}

\newcommand{\N}{\mathcal{N}}
\newcommand{\E}{\mathcal{E}}

\renewcommand{\footnoterule}{%
  \kern -3pt
  \hrule width 0.49 \textwidth height 0.5pt
  \kern 1pt
}
\usepackage[absolute]{textpos}
\TPMargin{6pt}

\begin{document}
\title{Online Distributed Optimization in Radial Power Distribution Systems: Closed-Form Expressions\\

}
\vspace{-3pt}
\author{\IEEEauthorblockN{Rabayet Sadnan$^{*1}$, Tom Asaki$^2$ and Anamika Dubey$^1$}\\
\vspace{-10pt}
\IEEEauthorblockA{{Department of Electrical Engineering and Computer Science$^1$, Mathematics and Statistics$^2$} \\
{Washington State University} \\
Email: \{rabayet.sadnan$^*$, tasaki, anamika.dubey\}@wsu.edu}
\vspace{-11pt}}

\maketitle
\begin{abstract}
The limitations of centralized optimization methods in managing power distribution systems operations motivate distributed control and optimization algorithms. However, the existing distributed optimization algorithms are inefficient in managing fast varying phenomena, resulting from highly variable distributed energy resources (DERs). Related online distributed control methods are equally limited in their applications. They require thousands of time-steps to track the network-level optimal solutions, resulting in slow performance. We have previously developed an online distributed controller that leverages the system's radial topology to achieve network-level optimal solutions within a few time steps. However, it requires solving a node-level nonlinear programming problem at each time step. This paper analyzes the solution space for the node-level optimization problem and derives the analytical closed-form solutions for the decision variables. The theoretical analysis of the node-level optimization problem and obtained closed-form optimal solutions eliminate the need for embedded optimization solvers at each distributed agent and significantly reduce the computational time and optimization costs.
\end{abstract}


\section{Introduction}
The development of optimal power flow methods (OPF) for electric power distribution systems has gained significant attention due to the increased number of distributed energy resources (DERs) with controllable smart inverters as they unfold the opportunity to operate power distribution grids more efficiently \cite{momoh1999review1,castillo2013survey}. To this end, both centralized and distributed OPF computational paradigms have emerged as viable mechanisms for distribution systems. The computational challenges posed by the centralized optimization paradigm and its susceptibility to a single-point failure motivate the use case for the distributed OPF algorithms. However, as the state-of-the-art distributed OPF methods are typically open-loop and slow in response, the intermittent nature of the DER generation makes it challenging to track the network-level optimal solutions. The primary drawback remains the requirement of a large number of message-passing rounds among the agents (on the order of $10^2 - 10^3$) to converge for a single-step optimization \cite{erseghe2014distributed, dall2013distributed,millar2016smart,magnusson2015distributed}. Note that a large number of communication rounds/message-passing events among distributed agents is not preferred since this leads to significant delays in decision-making. 


To alleviate these challenges, recent work proposes distributed online feedback-based voltage controllers to solve the OPF problem in a distributed manner. \cite{bolognani2014distributed, cavraro2017local,bernstein2019real,bastianello2020distributed,qu2019optimal,hu2019branch}. In contrast with traditional distributed optimization methods, these controllers do not wait to optimize for a time-step but asymptotically arrive at the optimal solution over several steps of real-time decision-making. They generally take one step towards the optimal solution and then move on to the next time step of the system simulations/observation. For example, in  \cite{bolognani2014distributed,cavraro2017local}, the proposed algorithm minimizes the active power loss upon implementing one step of the gradient descent method at each controllable node. In \cite{bernstein2019real}, a bilevel real-time controller is proposed that enables the agents to pursue a given performance objective and maintain operational limits by using a primal-dual projected gradient method. In \cite{bastianello2020distributed}, authors developed a proximal gradient method for online convex optimization, and in \cite{qu2019optimal}, a distributed controller is proposed that can limit the voltage, satisfy $Q$ capacity, and minimize a cost function. 

Although the existing real-time feedback-based online control algorithms can manage fast-changing system conditions, they pose several limitations. First, they are unable to reach global optimal solutions \cite{bolognani2014distributed,cavraro2017local}. Second, some of the existing methods require a central coordinator to update the global variable or to start an update sequence \cite{bastianello2020distributed, hu2019branch}; thus, they are not fully distributed. Moreover, these techniques require several time-steps to reach the optimal solution for a system with steady systems parameters \cite{qu2019optimal}. Thus, they are slow at tracking the network-level optimal solutions and show suboptimal performance for fast-changing system conditions (such as DER generation variability). Also, since the intermediate iterates are not optimal, they generally violate critical system operating constraints. To mitigate some of these concerns, previously, we developed an online distributed voltage controller for radial distribution systems, ENDiCO controller, based on the equivalence of networks principle \cite{sadnan2020real,sadnan2021distributed}. The proposed approach leverages the radial topology of the power distribution system and the associated unique power flow properties to develop the distributed feedback-control algorithms. This controller reduces the number of time steps required to track the optimal solutions by order of magnitude. However, it requires solving a nonlinear programming problem at each time step, thus necessitating an embedded optimization solver at each controllable node.  

This paper aims to develop analytical solutions using closed-form expressions for the node-level optimization problems in ENDiCO controllers, namely ENDiCO CLF. This eliminates the need for embedded optimization solvers at each controllable node, replacing them with a relatively simple algebraic computation, and greatly reduces the computational time needed to obtain network-level optimal solutions. We detail the solution space for Volt-Var (VVC) and Volt-Watt control (VWC) problems and use that to obtain an analytical solution for the decision variables.
Note that the existing distributed voltage control algorithms typically solve one step of the network-level optimization problem in a distributed way. Contrary to the existing online distributed control algorithms, our approach reduces the variable space for the optimization problem at each node and obtains the optimal decisions using closed-form expressions. The simulation results validate that the proposed approach can successfully track the network-level optimal solutions while maintaining the distribution system's operating constraints.

\vspace{-3pt}
\section{System Modeling \& Problem Formulation}
In this paper $(\cdot)^{(t)}$ represents the variable at time step $t$; $|~.~|$ symbolizes cardinality of a discrete set, or the absolute value of a number;
$j = \sqrt{-1}$; $\underline{(.)}$ and $\overline{(.)}$ are used to denote the minimum and maximum limit of any quantity, respectively.
\subsection{Network \& DER Model}
We assume a radial single-phase power distribution network, where $\N$ and $\E$ denote the set of nodes and edges of the system. Here, edge $ij \in \E$ identifies the distribution lines connecting the ordered pair of buses $(i,j)$ and is weighted with the series impedance of the line -- represented by $r_{ij}+jx_{ij}$. The set of load buses and DER buses are denoted by $\N_L$ and $\N_{D}$, respectively.
Let, for a given node $j$, node $i$ be the unique parent node, and $\N_{jk}=[k_1,k_2,.., k_n]$ be the set of children nodes for node $j$.
We denote $v_j$ and $l_{ij}$ as the squared magnitude of voltage and current flow at node $j$ and in branch $\{ij\}$, respectively. The network is modeled using the nonlinear branch flow equations \cite{baran1989optimal2} shown in \eqref{eqModel_nonlin}. Here, $p_{L_j}+jq_{L_j}$ is the load connected at node $j$, $P_{ij}, Q_{ij} \in \R$ are the sending-end active and reactive power flows for the edge $ij$, and $p_{Dj}+jq_{Dj}$ is the power output of the DER connected at node $j \in \N_{D}$.

\vspace{-0.5 cm}
\begin{small}
\begin{IEEEeqnarray}{C C}
\small
\IEEEyesnumber\label{eqModel_nonlin} \IEEEyessubnumber*
P_{ij}-r_{ij}l_{ij}-p_{L_j}+p_{Dj}= \sum_{k:j \rightarrow k} P_{jk}   \label{eqModel_nonlin1}\\
Q_{ij}-x_{ij}l_{ij}-q_{L_j}+q_{Dj}= \sum_{k:j \rightarrow k} Q_{jk} \label{eqModel_nonlin2}\\
v_j=v_i-2(r_{ij}P_{ij}+x_{ij}Q_{ij})+(r_{ij}^2+x_{ij}^2)l_{ij}\label{eqModel_nonlin3}\\
v_il_{ij} = P_{ij}^2+Q_{ij}^2 \label{eqModel_nonlin4}
\end{IEEEeqnarray}
\end{small}
\vspace{-0.5 cm}


The DERs are modeled as Photovoltaic modules (PVs) interfaced using smart inverters, capable of two-quadrant operation. At any node $j \in \N_D$, for the Volt-Var Control (VVC) method, the real power generation by the DER, $p_{Dj}$ is assumed to be known (measured). The reactive power generation, $q_{Dj}$, is controllable and modeled as the decision variable. Let the rating of the DER connected at node $j \in \N_D$ be $S_{DRj}$, then the limits on $q_{Dj}$ are given by \eqref{DG_lim}.

\vspace{-0.03cm}
 \begin{small}
\begin{equation}
\IEEEyesnumber\label{DG_lim}
-\sqrt{S_{DRj}^2-p_{Dj}^2} \leq q_{Dj} \leq \sqrt{S_{DRj}^2-p_{Dj}^2}
\end{equation}
 \end{small}
On the contrary, for the Volt-Watt Control (VWC), $q_{Dj}$ is set to $0$, and $p_{Dj}$ is assumed to be controllable and can vary between $0$ and $S_{DRj}$, see \eqref{DG_lim2}.

\vspace{-0.2cm}
\begin{small}
\begin{equation}
\IEEEyesnumber\label{DG_lim2}
0 \leq p_{Dj} \leq S_{DRj}
\end{equation}
 \end{small}
\vspace{-0.6cm}

\subsection{Distributed Real-Time Controller}
\begin{figure}[t]
    \centering
    \includegraphics[width=0.35\textwidth]{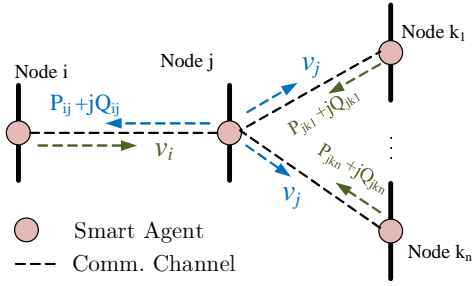}
        \vspace{-0.2cm}
    \caption{Cyber Layer Communication}
    \label{cyber_network}
\vspace{-0.8cm}
\end{figure}

Recently, we have developed a real-time, feedback-based distributed controller, ENDiCO, to solve network level optimal power flow problems \cite{sadnan2020real}. Briefly, each node $j \in \N_D$ solves the OPF problem defined by P1 \eqref{nonlin_OPF}. Either $p_{Dj}$ or $q_{Dj}$ is controlled to minimize some cost/objective function $\mathbf {f}$. At time-step $t$, node $j$ receives node voltage $v_i^{(t-1)}$ from the parent node $i$, and power flows $P_{jk}^{(t-1)}+jQ_{jk}^{(t-1)}$ from all of the children nodes in set $\N_{jk}$ (Fig. \ref{cyber_network}). The ENDiCO controller assumes the parent node voltage and the power flow to the children node to be constant, and solves the problem P1 locally for the reduced network.
Note that the controller at node $j$ only requires the upstream node voltage and downstream active and reactive power flows for optimization.

\textbf{\textit{Assumption 1:}} All the nodes in the network have an agent that can measure its local power flow quantities (node voltages and line flows) and communicate with neighboring nodes.

\vspace{-0.4cm} 
\begin{small}
\begin{IEEEeqnarray}{C C}
\IEEEyesnumber\label{nonlin_OPF} \IEEEyessubnumber*
\text{\textbf{(P1)}}\hspace{0.4cm} \min \hspace{0.2cm} \mathbf {f}^{(t)} \hspace{0.6cm}\\
P_{ij}^{(t)}-r_{ij}l_{ij}^{(t)}-p_{L_j}^{(t)}+p_{Dj}^{(t)}= \sum_{k:j \rightarrow k} P_{jk}^{(t-1)}  \label{eqModel_nonlin1}\\
Q_{ij}^{(t)}-x_{ij}l_{ij}^{(t)}-q_{L_j}^{(t)}+q_{Dj}^{(t)}= \sum_{k:j \rightarrow k} Q_{jk}^{(t-1)} \label{eqModel_nonlin2}\\
v_j^{(t)}=v_i^{(t-1)}-2(r_{ij}P_{ij}^{(t)}+x_{ij}Q_{ij}^{(t)})+(r_{ij}^2+x_{ij}^2)l_{ij}^{(t)}\hspace{0.4cm}\label{eqModel_nonlin3}\\
l_{ij}^{(t)} = \frac{(P_{ij}^{(t)})^2+(Q_{ij}^{(t)})^2}{v_i^{(t-1)}} \label{eqModel_nonlin4}\\
\underline{v} \leq v_j^{(t)} \leq \overline{v} \label{nonlin_OPF1}\\
l_{ij}^{(t)} \leq \left(I^{rated}_{ij}\right)^2  \label{nonlin_OPF2}\\
\text{DG Limit:} \hspace{0.3 cm}\text{equation}\hspace{0.1cm}\eqref{DG_lim} \hspace{0.1 cm} \text{or} \hspace{0.1 cm} \eqref{DG_lim2}\label{lim_DG}
\end{IEEEeqnarray}
\end{small}
\vspace{-0.4cm}

Here, $\underline{v} = 0.95^2$ and $\overline{v} = 1.05^2$ pu are the limits on bus voltages, and $I^{rated}_{ij}$ is the thermal limit for the branch $\{ij\}$. The convergence of the boundary variables is obtained using the Fixed-Point Iteration (FPI) update denoted by \eqref{fpi_up}. Here, the variable $X$ can be the power flow requirements from the child nodes, i.e., $\{P_{jk},Q_{jk}\}$, or the voltage at the parent node, $v_i$. Further, instead of a constant value, $\alpha$ can be made adaptive.

\vspace{-0.2cm}
\begin{small}
\begin{equation}\label{fpi_up} 
X^{(t-1)}:= \frac{X^{(t-1)}+\alpha X^{(t-2)}}{1+\alpha}
\end{equation}
\end{small}
\vspace{-0.3cm}

Let $\mathbf {F}^{(t)}$ be the constrained optimization problem defined by P1 at time step $t$. For VVC, node $j$ attains the optimal reactive power dispatch $q_{Dj}^{(t)}$ using \eqref{dist_control1}. For the VWC, the optimal active power dispatch $p_{Dj}^{(t)}$ is solved using \eqref{dist_control2}. The resulting nonlinear optimization problem $\mathbf {F}^{(t)}$ is in five variables; $\{P_{ij}^{(t)}, Q_{ij}^{(t)}, v_{j}^{(t)}, l_{ij}^{(t)}, q_{Dj}^{(t)}\}$ for VVC, and $\{P_{ij}^{(t)}, Q_{ij}^{(t)}, v_{j}^{(t)}, l_{ij}^{(t)}, p_{Dj}^{(t)}\}$ for VWC. All controllable nodes in the system solve P1 in parallel and calculate their respective dispatches using \eqref{dist_control}. Further details can be found in \cite{sadnan2020real}.

\vspace{-0.3cm}


\begin{small}
\begin{IEEEeqnarray}{C C}
\IEEEyesnumber\label{dist_control} \IEEEyessubnumber*
\text{For VVC: }\hspace{0.2cm} q_{Dj}^{(t)} = q_{Dj}^* =  \underset{q_{Dj}}{\arg\min} \hspace{0.2cm} \mathbf {F}^{(t)}(q_{Dj})\label{dist_control1}\\
\text{For VWC: }\hspace{0.2cm} p_{Dj}^{(t)} = p_{Dj}^* =  \underset{p_{Dj}}{\arg\min} \hspace{0.2cm} \mathbf {F}^{(t)}(p_{Dj}) \label{dist_control2}
\end{IEEEeqnarray}
\end{small}
\vspace{-0.5cm}

\section{Distributed Controller with Closed-Form Solution: ENDiCO CLF}
In this section, we obtain the closed-form solutions for the distributed OPF problem detailed in P1 (4).
The controllable nodes calculate the optimal decision variables for the current time-step using the obtained closed-form expressions in this section. First, we detail the solution space and the closed-form solution for the VVC case. Next, the analytical solution for the VWC problem is constructed. The ENDiCO CLF algorithm is discussed last.\\
\textbf{\textit{Assumption 2:}} The loads are modeled as constant power loads. 

\subsection{Analytical Solution for Volt-Var Control} For the VVC, let us define the optimization variables of the optimization problem P1 as $x = \{x_1,x_2,x_3,x_4,x_5\} = \{P_{ij}^{(t)}, Q_{ij}^{(t)}$, $ v_{j}^{(t)}, l_{ij}^{(t)}, q_{Dj}^{(t)}\}$. We denote the constants (at time-step $(t-1))$ as the following: $P = \sum_{k:j \rightarrow k} P_{jk}^{(t-1)} +p_{L_j}^{(t)} -p_{Dj}^{(t)}$, $Q = \sum_{k:j \rightarrow k} Q_{jk}^{(t-1)} +q_{L_j}^{(t)}$, $V = v_i^{(t-1)}$  and $z = z_1+jz_2 = r_{ij}+jx_{ij}$. Then, we can write the aforementioned OPF problem in P1 using \eqref{gen_opf}. Here, $l =\{l_1,l_2,l_3,l_4,l_5\}$ and $u =\{u_1,u_2,u_3,u_4,u_5\}$ are the lower and upper bounds for each of the problem variables.

\vspace{-0.4cm}
\begin{small}
\begin{IEEEeqnarray}{C C}
\IEEEyesnumber\label{gen_opf} \IEEEyessubnumber*
\text{\textbf{(P2)}} \hspace{0.4cm} \min \hspace{0.2cm} \mathbf{f}(x) \hspace{0.6cm}\\
x_1 = z_1x_4+P    \label{gen_nonlin1}\\
x_2= z_2x_4+Q-x_5  \label{gen_nonlin2}\\
x_3 =V-2(z_1x_1 +z_2x_2 )+z^2x_4 \hspace{0.4cm}\label{gen_nonlin3}\\
Vx_4 = x_1^2+x_2^2 \label{gen_nonlin4}\\
l\leq x \leq u \label{gen_nonlin5}
\end{IEEEeqnarray}
\end{small}
\vspace{-0.4cm}

The general solution for the set of linear equalities \eqref{gen_nonlin1}-\eqref{gen_nonlin3} is represented in parametric form as follows.

\vspace{-0.2cm}
\begin{equation}
\begin{bmatrix} x_1\\
x_2\\x_3\\x_4\\x_5\end{bmatrix}= 
\begin{bmatrix} P\\
Q\\V-2(z_1P +z_2Q)\\0\\0\end{bmatrix}+
\begin{bmatrix} z_1\\
z_2\\-z^2\\1\\0\end{bmatrix}x_4+
\begin{bmatrix} 0\\
-1\\2z_2\\0\\1\end{bmatrix}x_5
\label{lin_sol}
\end{equation}

Using $x_1$ and $x_2$ from \eqref{lin_sol} in the non-linear constraint \eqref{gen_nonlin4}, we get \eqref{elp1} which represents the solution space of the local problem. This can further be simplified to \eqref{elp2}.

\vspace{-0.3cm}
\begin{small}
\begin{IEEEeqnarray}{C C}
\IEEEyesnumber \IEEEyessubnumber*
(P+z_1x_4)^2+(Q+z_2x_4-x_5)^2 = Vx_4    \label{elp1}\\
\begin{split}
\Rightarrow z^2x_4^2-2z_2x_4x_5+x_5^2+x_4(2Pz_1+&2Qz_2-V)\\+x_5(-2Q)+P^2+Q^2 &= 0
\end{split}
\label{elp2}
\end{IEEEeqnarray}
\end{small}
\vspace{-0.3cm}

Equation \eqref{elp2} represents an ellipse whose general form is $Ax_4^2+Bx_4x_5+Cx_5^2+Dx_4+Ex_5+F = 0$, where, $A = z^2$, $B = -2z_2$, $C =1$, $D=2Pz_1+2Qz_2-V $, $E=-2Q$, and $F=P^2+Q^2$. The angle $\theta$ (w.r.t. the $x_4$ axis) and the center $(x_4^0,x_5^0)$ of the ellipse can be found using \eqref{theta} and \eqref{el_cent}.

\vspace{-0.4cm}
\begin{small}
\begin{eqnarray}\label{theta}
\theta &=& \hspace{0.2cm} \text{arctan}\hspace{0.2 cm} \left(-\frac{1}{2z_2}(1-z^2-\sqrt{(z^2-1)+4z_2^2})\right)
\end{eqnarray}
\vspace{-0.3cm}
\begin{IEEEeqnarray}{C C}
\IEEEyesnumber\label{el_cent} \IEEEyessubnumber*
x_4^0 = \frac{2CD-BE}{B^2-4AC} =\frac{V-2Pz_1}{2z_1^2}\\
x_5^0 = \frac{2AE-BD}{B^2-4AC} =\frac{Vz_2+2Qz_1^2-2Pz_1z_2}{2z_1^2}
\end{IEEEeqnarray}
\end{small}
\vspace{-0.3cm}

\begin{figure}[t]
\centering
\subfloat[\centering{VVC Case}]{\includegraphics[width=0.23\textwidth]{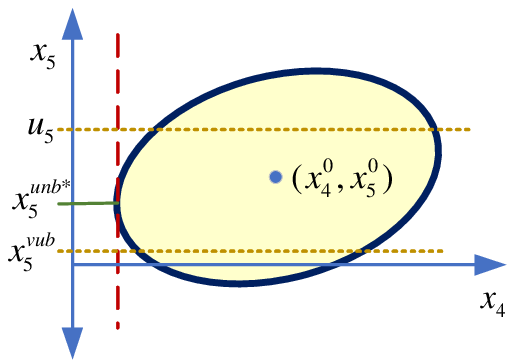}}
\subfloat[\centering{VWC Case}]{\includegraphics[width=0.23\textwidth]{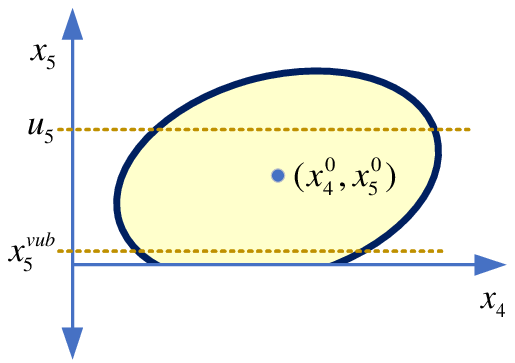}}
\caption{Solution Space of the Distributed OPF}
\label{ellipse_12}
\vspace{-0.7cm}
\end{figure}

Now that we have the solution space ready for VVC (Fig. \ref{ellipse_12}a), we discuss a specific VVC problem where we control reactive power generation of DERs, $x_5$ to minimize the active power losses, $r_{ij}l_{ij}$. This is equivalent to minimizing $\mathbf{f}(x) = x_4$. Upon removing the box constraints of the variables, we first attain the optimum solution for the relaxed problem. Then the solution is projected onto the feasible space defined by the bounds of the variable space. The relaxed solution, $x_5^{unb*}$, is obtained upon differentiating the left-hand side of \eqref{elp1} w.r.t. $x_5$ and setting the result to $0$ \eqref{min_x5}. After substituting the $x_5^{unb*}$ in \eqref{elp1}, we also get the minimum current flow, $x_4^{unb*}$, for the relaxed problem using (14).

\vspace{-0.4cm}
\begin{small}
\begin{eqnarray}
\frac{d}{dx_5}\bigg\{(P+z_1x_4)^2+(Q+z_2x_4-&x_5)^2\bigg\} = 0\label{dif_x5}\\
\Rightarrow  x_5^{unb*} = Q+z_2x_4^{unb*}&\label{min_x5}\\
x_4^{unb*} = \frac{V-2Pz_1-\sqrt{V^2-4VPz_1}}{2z_1^2}
\label{min_x4}
\end{eqnarray}
\end{small}
\vspace{-0.4cm}

Next, $x_5^{unb*}$ is projected onto the feasible space by applying the voltage bounds and DERs physical limits. For active power loss minimization problem, the voltage will have an active upper bound, i.e., $x_3\le {u_3}$. Thus the projected value, $x_5$, onto the feasible space for voltage bounds, $x_5^{vub}$, can be found using \eqref{vub}. Here, $x_4$ is approximated as, $x_4\approx \frac{P^2+(Q-x_5)^2}{V}$.

\vspace{-0.3cm}
\begin{small}
\begin{equation}
\begin{split}
&V-2(z_1P+z_2Q)-z^2x_4+2z_2x_5\le u_3\\
\Rightarrow &~\frac{z^2}{2z_2V}(x_5)^2 - \left(\frac{z^2Q}{z_2V}+1\right)x_5+Q+\frac{z_1}{z_2}P+\\&\frac{1}{2z_2}\left(u_3-V\right) +\frac{z^2}{2z_2V}\left(P^2+Q^2\right) \ge 0 
\\\Rightarrow &~x_5 \le \frac{-b_1-\sqrt{b_1^2-4a_1c_1}}{2a_1} = x_5^{vub}
\end{split}
\label{vub}
\end{equation}
\end{small}
\vspace{-0.1cm}
\noindent where, $a_1 = {z^2/}{2z_{2}V}$, $b_1= -({z^2Q/}{z_2V}+1)$ and $c_1 = Q+\frac{z_1}{z_2}P+\frac{1}{2z_2}(u_3-V)+\frac{z^2}{2z_2V}(P^2+Q^2) $. The other solution would be $x_5\ge \frac{-b_1+\sqrt{b_1^2-4a_1c_1}}{2a_1}$, which is at the far right of the ellipse, and thus not a feasible solution. After considering the physical limits of the DER reactive power generation ($u_5$), the analytical solution for the power loss minimization problem is given by \eqref{x_5st}. 
A similar analysis can be performed for other VVC problems to obtain a closed-form solution. One example could be minimizing bus voltages, $x_3$. This problem will activate lower bound constraints on the voltage ($x_3\ge l_3$), instead of the upper bound constraint.

\vspace{-0.1cm}
\begin{small}
\begin{equation}
q_{Dj}^{(t)} = x_5^* = \text{min}\hspace{0.2 cm}\{x_5^{unb*},x_5^{vub},u_{5}\}
\label{x_5st}
\end{equation}
\end{small}
\vspace{-0.6cm}

\subsection{Analytical Solution for Volt-Watt Control}
Similarly for the VWC case, we use the same constants $Q,V,z$ and same variables $x_1$ through $x_4$. However, we define $P = \sum_{k:j \rightarrow k} P_{jk}^{(t-1)} +p_{L_j}^{(t)}$ and $x_5 = p_{Dj}^{(t)}$. The optimization problem is formulated using \eqref{gen_opf_vwc}

\vspace{-0.4cm}
\begin{small}
\begin{IEEEeqnarray}{C C}
\IEEEyesnumber\label{gen_opf_vwc} \IEEEyessubnumber*
\text{\textbf{(P3)}} \hspace{0.4cm} \min \hspace{0.2cm} \mathbf{f}(x) \hspace{0.6cm}\\
x_1 = z_1x_4+P-x_5 \label{gen_nonlin6}\\
x_2= z_2x_4+Q  \label{gen_nonlin7}\\
\text{constraint}\hspace{0.2cm} \eqref{gen_nonlin3},\eqref{gen_nonlin4},\eqref{gen_nonlin5}
\end{IEEEeqnarray}
\end{small}
\vspace{-0.5cm}

Using similar methods, we obtain the solution space (the ellipse in Fig. \ref{ellipse_12}b) for the VWC problem as well.
Next, we detail the closed-form solution for a specific VWC problem objective. Specifically, we detail the closed-form solution for the problem of minimizing DER curtailment or maximizing the DER power generation. Thus, the cost function (for individual node) is given by, $\mathbf{f}(x) = -x_5$. Since we want to maximize $x_5$, the limiting constraint will be either the upper bound of the voltage $u_3$, or the physical limit of the DER power generation, $u_5$. By activating the constraint $x_3\le u_3$, we get $x_5^{vub}$-- the maximum amount of DER generation that does not violate the voltage bounds. A similar method detailed in \eqref{vub} is used to obtain, $x_5^{vub}$, expressed using \eqref{vub2}.

\vspace{-0.25cm}
\begin{figure}[b]
\vspace{-0.25cm}
    \centering
    \includegraphics[width=0.4\textwidth]{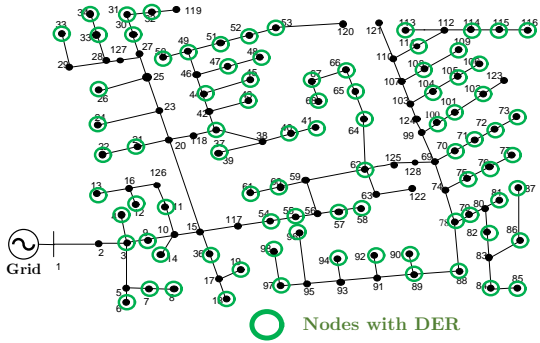}
        \vspace{-0.2cm}
    \caption{IEEE 123 Bus Test System with 85 DERs}
    \label{123_bus}
\end{figure}

\begin{small}
\begin{equation}\label{vub2}
x_5\le \frac{-b_2-\sqrt{b_2^2-4a_2c_2}}{2a_2} = x_5^{vub}
\end{equation}
\end{small}
\vspace{-0.3cm}

\begin{figure*}[t]
\centering
\hspace{-0.3cm}
\subfloat[\centering{Normalized Load and PV generation}]{\includegraphics[width=0.32\textwidth]{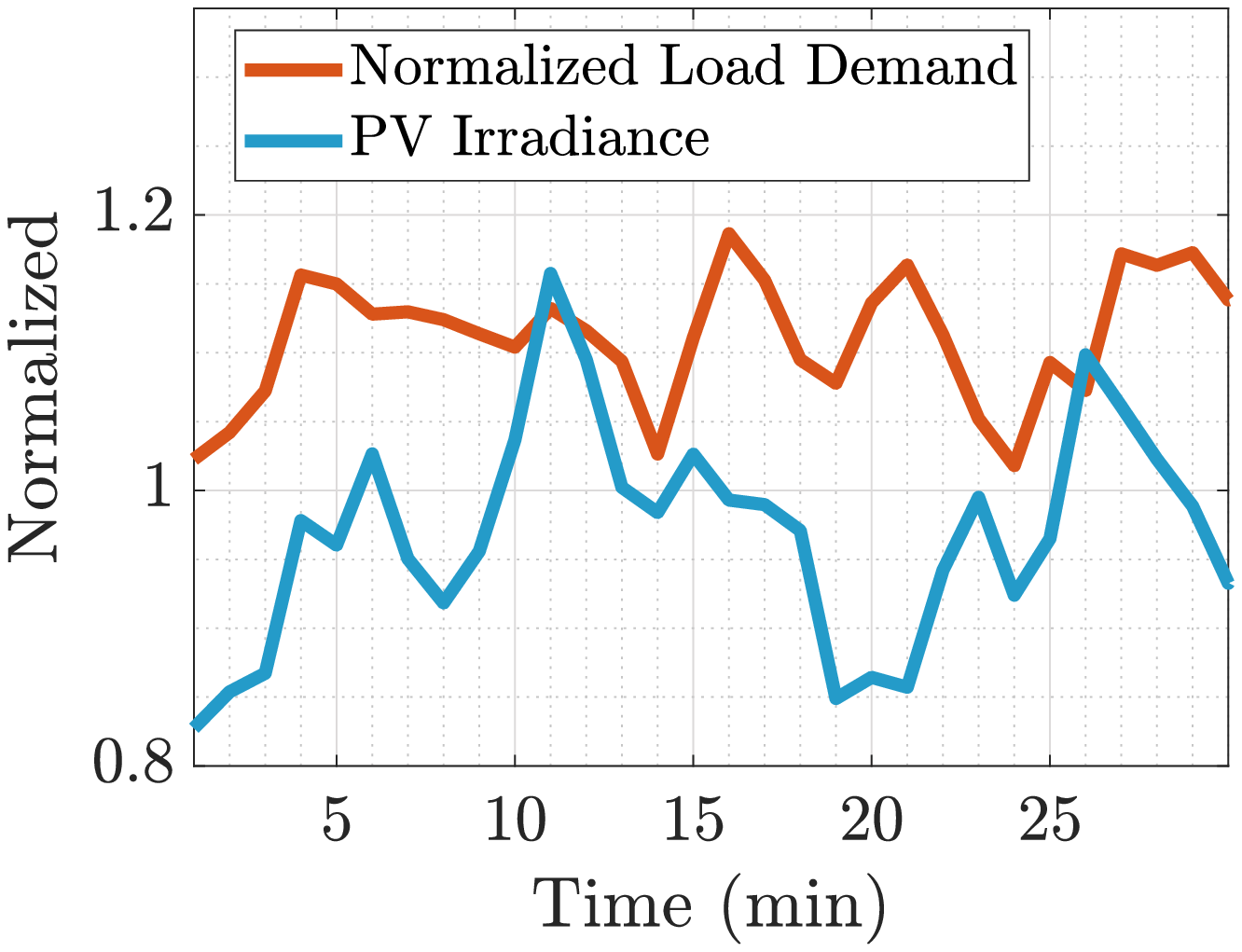}\label{PV_var}}
\hspace{-0.4cm}
\subfloat[\centering{VVC: Minimized Active Power Loss}]{\includegraphics[width=0.32\textwidth]{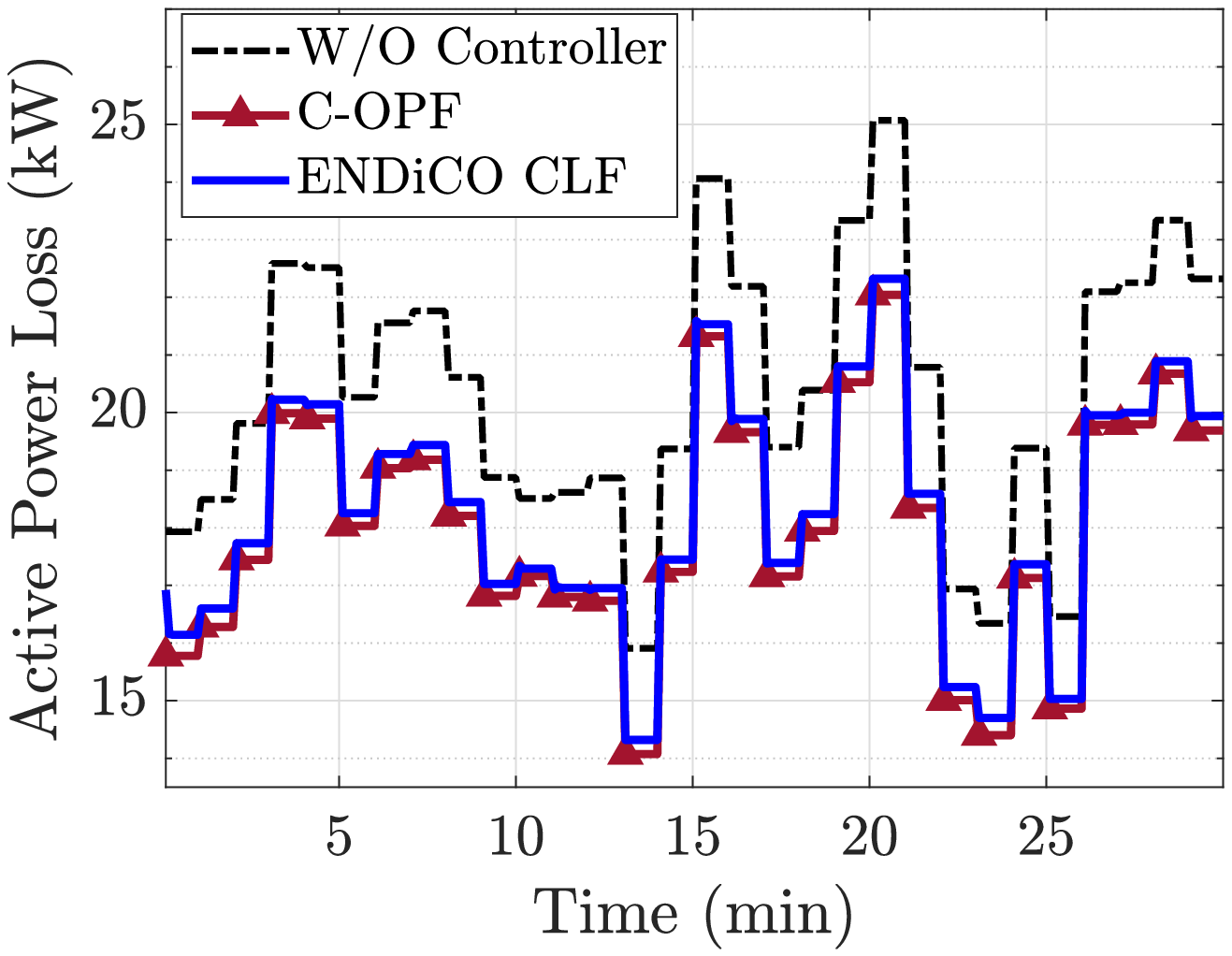}
\label{Loss_track}}
\subfloat[\centering{VWC: Maximized DER Generation}]{\includegraphics[width=0.32\textwidth]{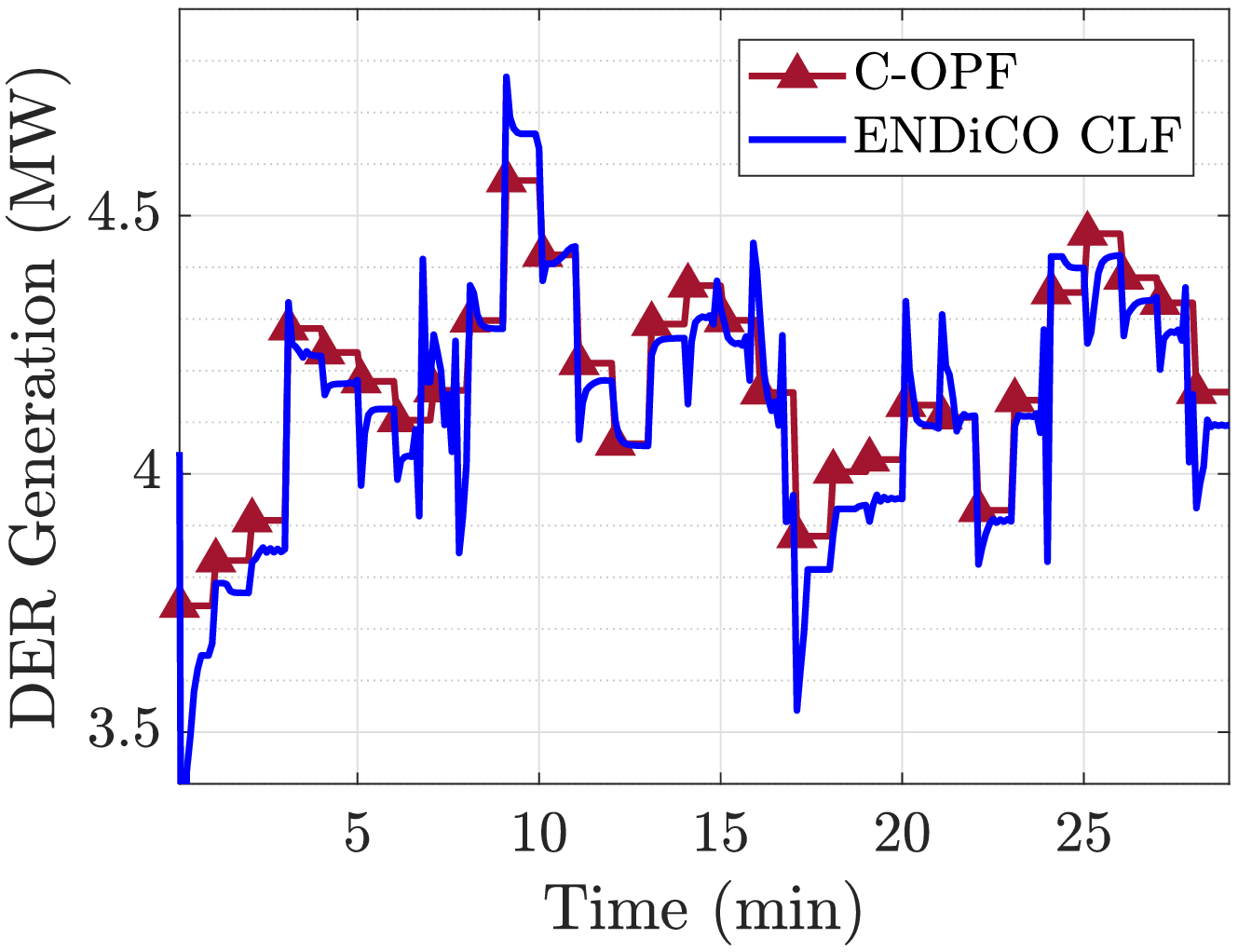}\label{Loss_track2}}
    \vspace{-0.2cm}
\caption{Time-Series Result for ENDiCO CLF}
\label{Result_obj}
\vspace{-0.4cm}
\end{figure*}
\noindent where, $a_2 = \frac{z^2}{2z_1V}$, $b_2= -\left(\frac{z^2P}{z_1V}+1\right)$ and $c_2 = P+\frac{z_2}{z_1}Q+\frac{1}{2z_1}(u_3-V)+\frac{z^2}{2z_1V}(P^2+Q^2)$. After considering the physical limits on DER active power generation, $x_5\le u_5$, we obtain the closed-form expression for the optimal active power generation from DER as \eqref{x_5stvw}.

\vspace{-0.2cm}
\begin{small}
\begin{equation}
p_{Dj}^{(t)} = x_5^* = \text{min}\hspace{0.2 cm}\{x_5^{vub},u_{5}\}
\label{x_5stvw}
\end{equation}
\end{small}
\vspace{-0.6cm}

{\color{black}
\subsection{Online Feedback-based ENDiCO-CLF Controller}
}
In this section, the proposed ENDiCO CLF algorithm is discussed. This method is based on the previously developed ENDiCO controller, and we use the same communication infrastructure; however, we incorporate the analytical solutions for finding the local optimal set-points. The algorithm of the ENDiCO CLF controller is detailed in Algorithm 1.

\vspace{-6pt}
\begin{algorithm}[!h]
\small
\caption{\small ENDiCO CLF}\label{alg:PSLP}
\SetAlgoLined
\SetKwInOut{Rx}{Receive}
\SetKwInOut{Tx}{Transmit}
\SetKwInOut{St}{Calculate}
\SetKwInOut{ND}{Node}
\SetKwInOut{TE}{Time Step }
\SetKwInOut{Stp}{Steps }
\ND {$\forall  j \in \N_D$}
\TE{t}
\Rx{$v_i^{(t-1)}$ and $P^{(t-1)}_{jk_i}+jQ^{(t-1)}_{jk_i}$}
\Tx{$v_j^{(t)}$ and $P^{(t)}_{ij}+jQ^{(t)}_{ij}$ }
\vspace{3pt}
\Stp{ }
\vspace{2pt}
\vspace{2pt}
Calculate $\sum P_{jk}^{(t-1)}+jQ_{jk}^{(t-1)}$ from all the $P^{(t-1)}_{jk_i}+jQ^{(t-1)}_{jk_i}$, received from child nodes $k_i \in \N_{jk}$ \\
\vspace{2pt}
Approximate the upstream and downstream network of line $\{ij\}$ with fixed value of $v_i^{(t-1)}$ and $\sum P_{jk}^{(t-1)}+jQ_{jk}^{(t-1)}$\\
\vspace{3pt}
Solve optimization problem \eqref{dist_control} using closed-form solutions; equation \eqref{x_5st} \& \eqref{x_5stvw} for VVC \& VWC, respectively
\\
\vspace{2pt}
Implement the set point $q^{(t)}_{Dj}$ for VVC, and $p^{(t)}_{Dj}$ for VWC at node $j$\\
\vspace{2pt}
Measure the node voltage $v_j^{(t)}$ at node $j$ and complex power flow $P^{(t)}_{ij}+jQ^{(t)}_{ij}$ in the line $\{ij\}$\\
\vspace{2pt}
Send $v_j^{(t)}$ and $P^{(t)}_{ij}+jQ^{(t)}_{ij}$ to child nodes $k_i$ and parent node $i$, respectively\\
\vspace{2pt}
Receive $v_i^{(t)}$ and $P^{(t)}_{jk_i}+jQ^{(t)}_{jk_i}$ from parent and child nodes, respectively\\
\vspace{2pt}
Move forward to the next time step $(t+1)$
\vspace{4pt}
	\label{algo}
\end{algorithm}




\vspace{-0.4cm}
\section{Case Studies}
In this section, we validate the proposed ENDiCO CLF
for both VVC and VWC problems using several test cases.
The performance and the solution quality for the proposed controller is analyzed by evaluating (i) its ability to track the optimal solutions, (ii) nodal voltage control capability, 
and (iii) time required to achieve the optimal solutions. The proposed ENDiCO CLF controller is also compared against a centralized OPF (C-OPF) solution via simulations.
\vspace{-0.1cm}
\subsection{Test Systems}
\vspace{-0.1cm}
We select IEEE-123 bus test system with $100$\% customer PV penetration ($85$ DERs) for simulation studies (Fig. \ref{123_bus}). The normalized load curves and high PV irradiance for simulated $30$-minute test case is shown in Fig. \ref{PV_var}. It is assumed that the load and PV irradiance change at every $1$-minute time-interval.
The time resolution for the controller is assumed $6$ seconds, i.e., each controller takes an optimal action and communicates with its immediate neighbor within 6 sec -- this is in accordance with related literature in this domain \cite{magnusson2020distributed}.
This proposed ENDiCO CLF controller is evaluated for two types of OPF problems: (i) VVC: loss minimization
and (ii) VWC: DER curtailment.
The rated real power of the DERs are randomly chosen between $1 - 5kW$ and $12 - 60kW$ for VVC and VWC case, respectively. The maximum kVA rating for the PV panel is set at $120\%$ of their rated real power generation. Also, we assume $\alpha$ to be $0$ and $10$ for VVC and VWC case, respectively.



\vspace{-0.1cm}
\subsection{Simulation Results}
\vspace{-0.1cm}
The numerical simulation of the ENDiCO CLF controllers are detailed in this section. Here, we compare the objective function values and the nodal voltages. We also compare the solution times for different OPF methods.

\begin{figure}[t]
\vspace{-0.2cm}
\centering
\subfloat[\centering{VVC Case}]{\includegraphics[width=0.25\textwidth]{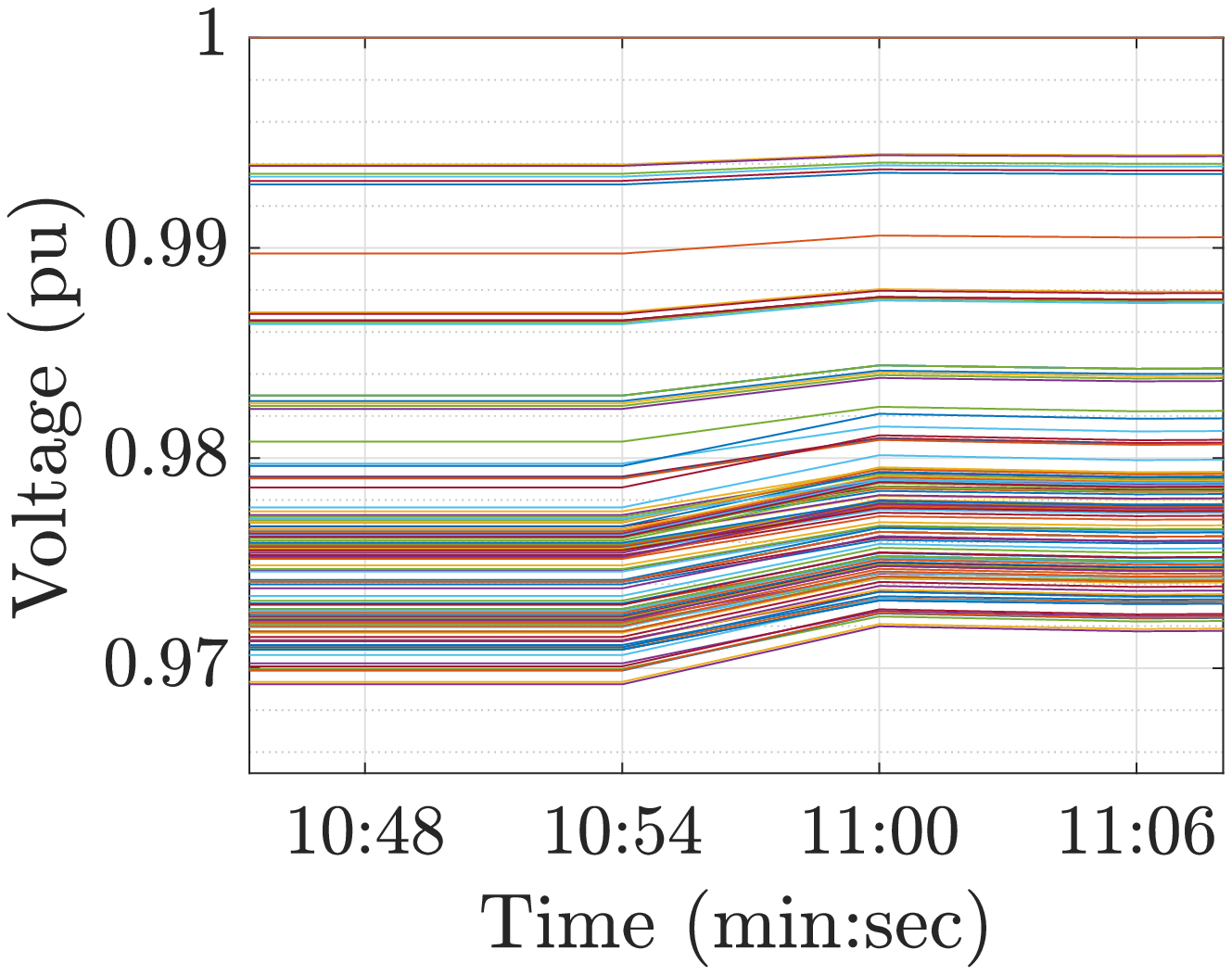}}
\subfloat[\centering{VWC Case}]{\includegraphics[width=0.25\textwidth]{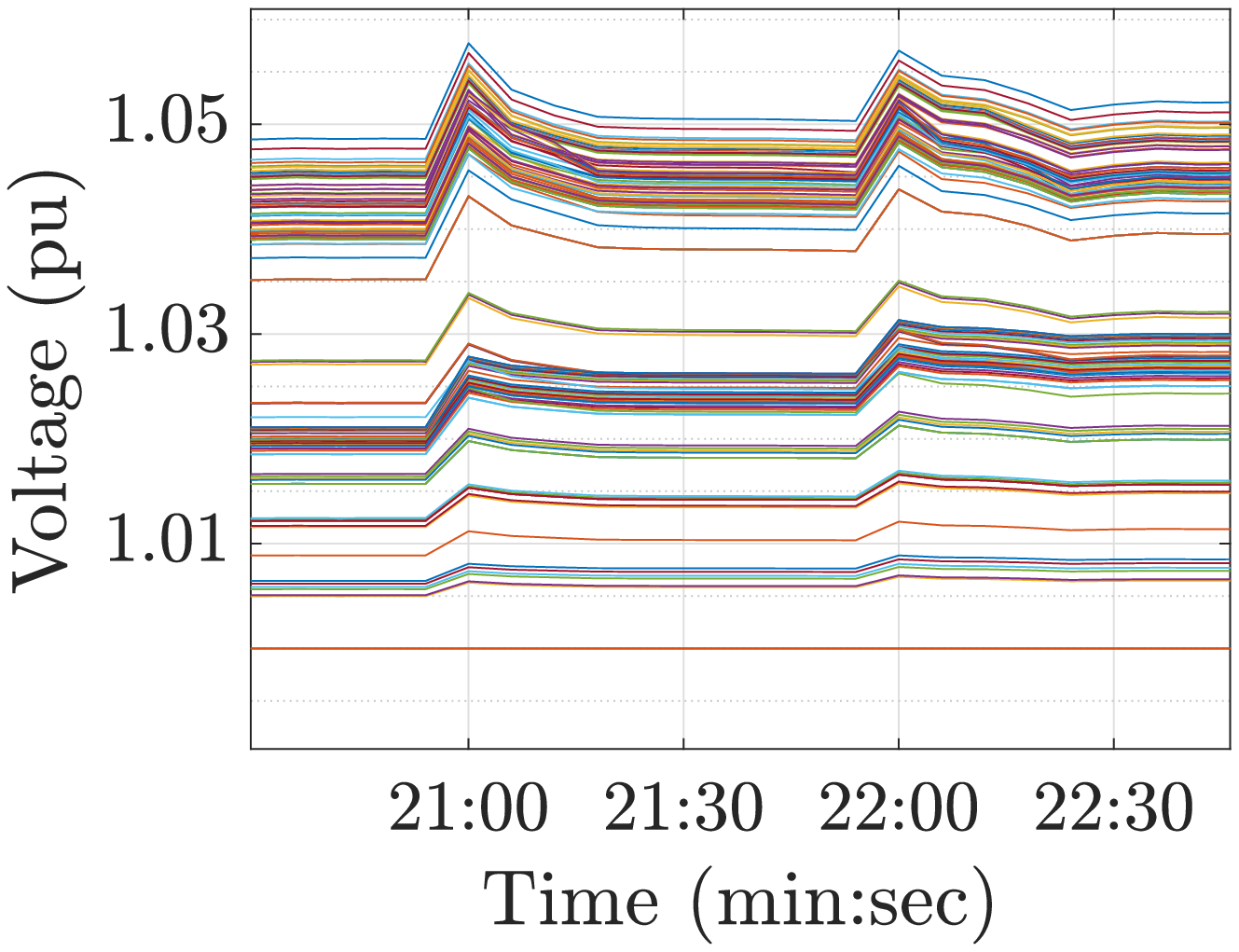}}
\caption{Voltage Oscillation for ENDiCO CLF}
\label{V_osc}
\vspace{-0.6cm}
\end{figure}

\subsubsection{Objective Value Tracking}
The Proposed ENDiCO CLF controller can track the value of the objective function for both VVC and VWC cases (Fig. \ref{Loss_track}, \ref{Loss_track2}). 
Here, C-OPF solves the same optimization problem centrally by collecting all the data from the network, and then dispatching the decision variables to the controllable nodes. From Fig. \ref{Result_obj}, it is evident that the developed ENDiCO controller can track the optimal solutions for the given time resolution.
Given these fast-changing scenarios, the proposed controller performs well and accurately tracks the original solution with minimum oscillations.

\begin{figure*}[t]
\centering
\hspace{-0.3cm}
\subfloat[\centering{VVC: ENDiCO CLF}]{\includegraphics[width=0.26\textwidth]{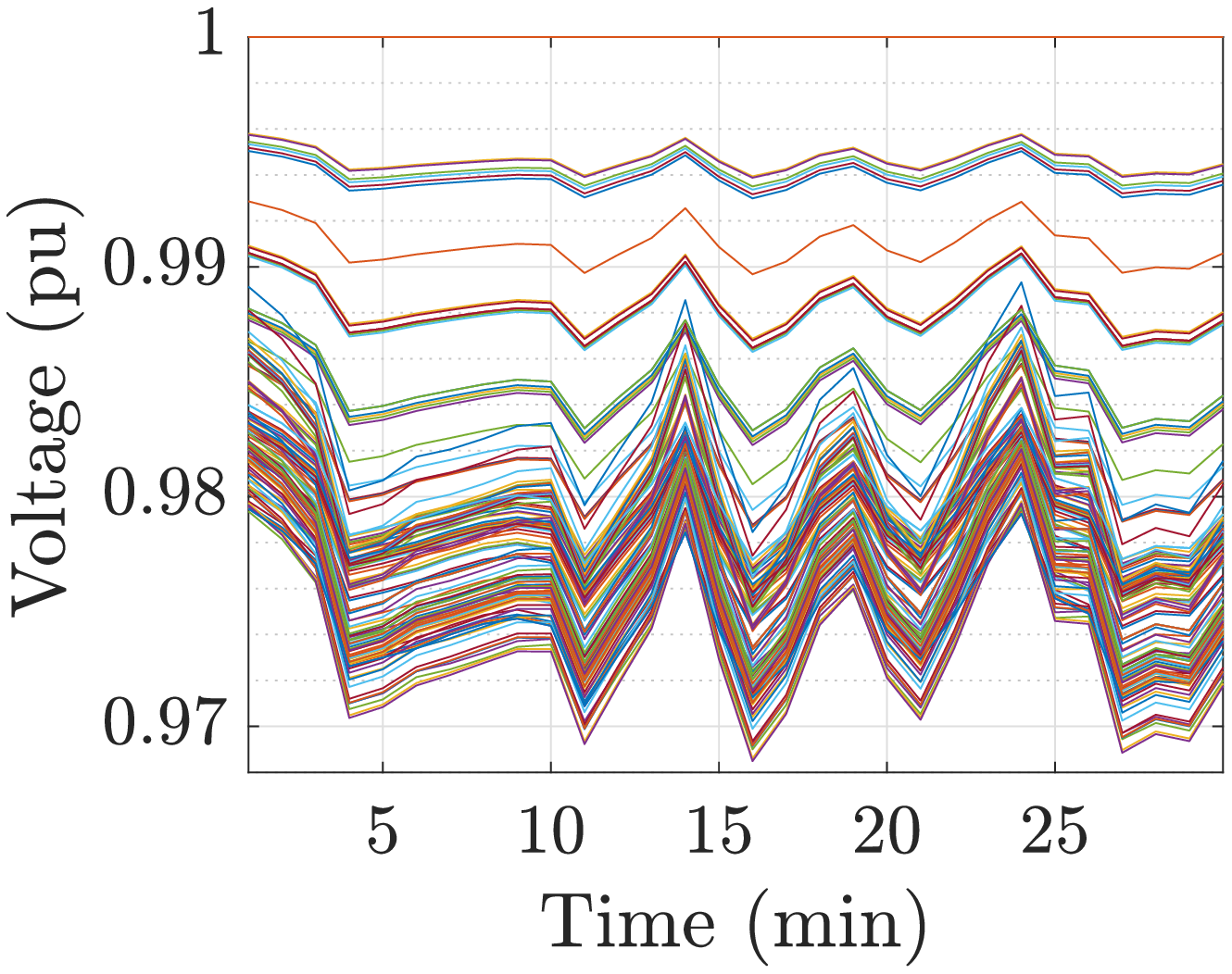}\label{LVC}}
\hspace{-0.4cm}
\subfloat[\centering{VVC: C-OPF}]{\includegraphics[width=0.26\textwidth]{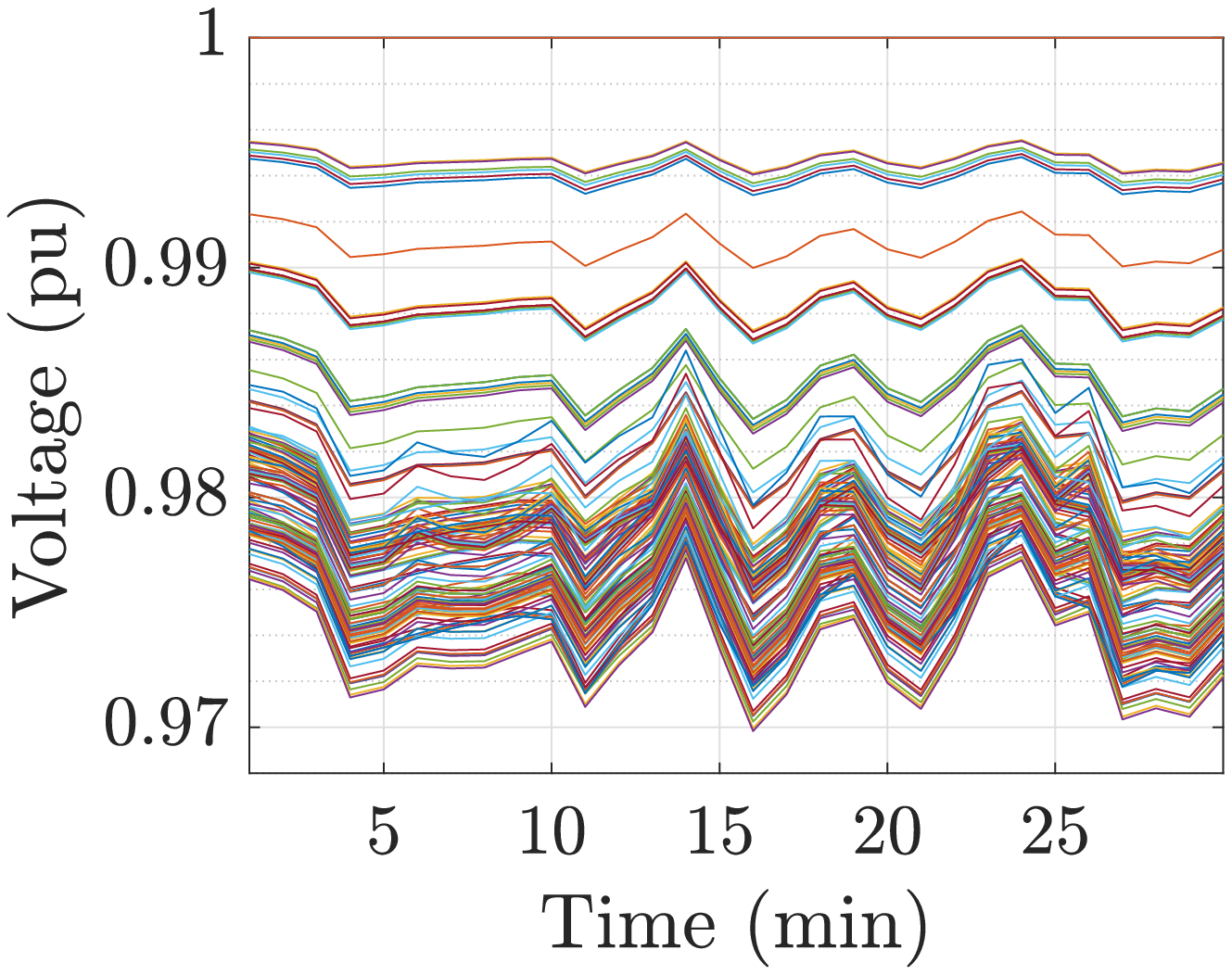}\label{LVNC}}
\hspace{-0.4cm}
\subfloat[\centering{VWC: ENDiCO CLF}]{\includegraphics[width=0.26\textwidth]{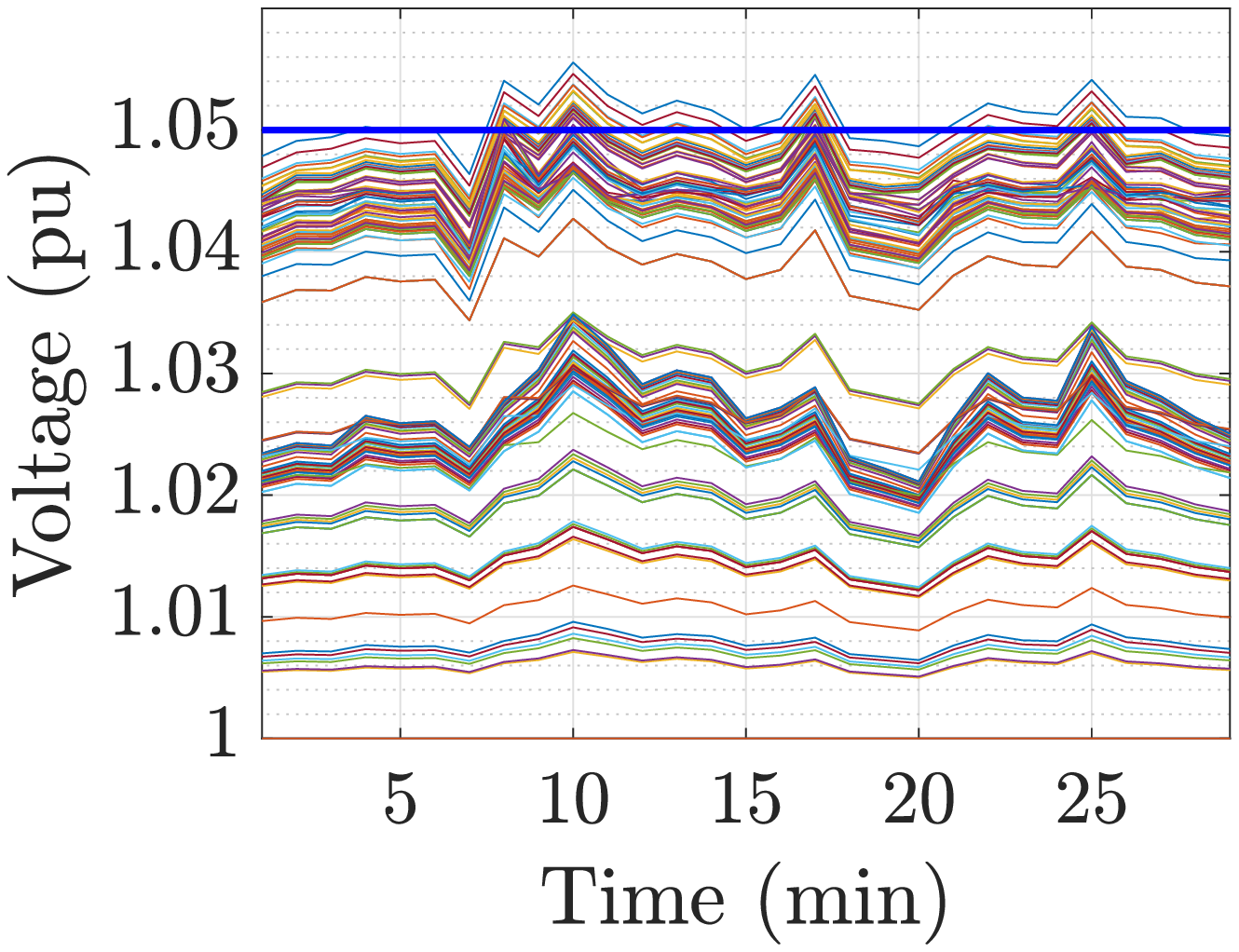}\label{PVC}}
\hspace{-0.4cm}
\subfloat[\centering{VWC: C-OPF}]{\includegraphics[width=0.26\textwidth]{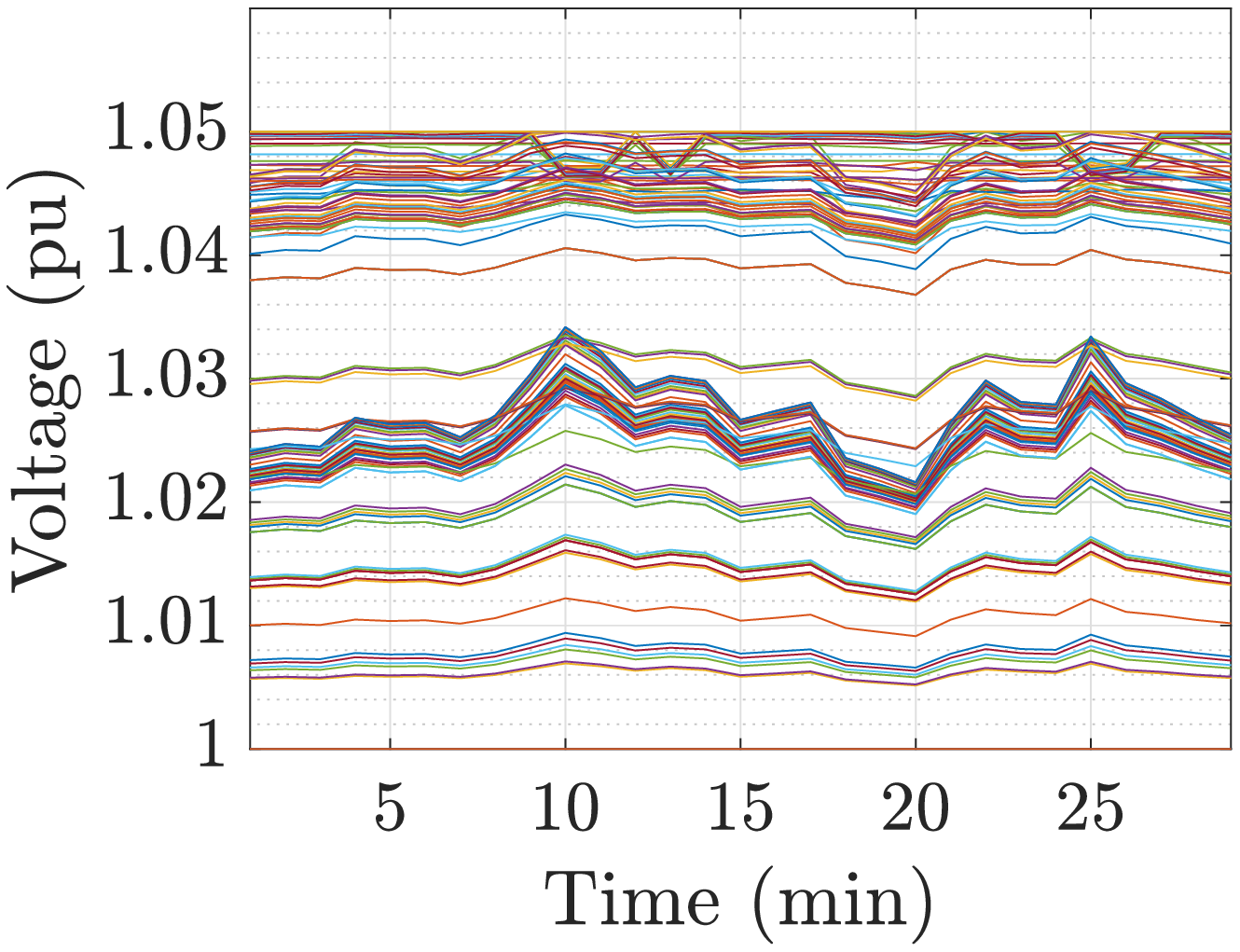}\label{PVNC}}
    \vspace{-0.2cm}
\caption{Comparison of nodal Voltage with Centralized OPF Solutions for Different Cases}
\label{Result_all}
\vspace{-0.4cm}
\end{figure*}

Fig. \ref{Loss_track} shows the active power loss in the line for the three cases: ENDiCO CLF, C-OPF and without any controller. Note that, for the VVC problem, the proposed ENDiCO CLF approach takes only $2$ time-steps to converge. 
The maximum tracking error in this case is not more than 0.3 kW. In Fig. \ref{Loss_track2}, the DER generation output for each minute is shown for the VWC problem.  Convergence is achieved within $4$ time-steps. An initial oscillation can be seen at the beginning of each time step, but is fully suppressed after two or three time steps. Also, the difference between the C-OPF solutions and the converged ENDiCO CLF solutions are less than $1.8\%$. It is to be noted that for the DER maximization problem (VWC), we have assumed an extreme DER deployment case where generation in the system is $\sim2.5 MW$ higher than the total load in the system. The proposed ENDiCO CLF performs reasonably well even during extreme scenarios. For more typical cases, the oscillations and tracking errors are significantly less.

\subsubsection{Nodal Voltage Comparison}
To evaluate the solution quality of ENDiCO CLF in terms of nodal voltages, we first discuss the voltage oscillations when system parameters change, and then we compare the stable voltage solution with C-OPF results.

\textit{(i) Voltage Oscillations:}
When the network parameters such as load and DER generation change, small oscillations (or no oscillations) can be seen in the node voltages (Fig. \ref{V_osc}). From the Fig. \ref{V_osc}a, we can see that for VVC with moderately stressed distribution system, the node voltage stabilizes fast after each change in the network parameters.
In addition to that, VWC with highly stressed system that has significantly high DER penetration, a small voltage oscillation is present at the very beginning of each time step (Fig. \ref{V_osc}b). The maximum oscillation recorded for the simulated case is of $0.005$ pu in the 30-min simulation window, but that reaches a stable value after $4$ time-steps. 

\textit{(ii) Voltage Control:}
Next the performance of ENDiCO CLF controller is evaluation for its ability to maintain the system voltages. Here, we compare the voltage profile obtain upon implementing the ENDiCO CLF controller with the results obtained from C-OPF. Fig. \ref{Result_all} shows the nodal voltages of ENDiCO CLF and C-OPF for both VVC and VWC cases. From the figure, it is clear that for VVC, nodal voltages for the ENDiCO CLF and C-OPF match closely. The maximum voltage difference between these two methods is only $0.001$ pu. Even for the highly stressed VWC case, we can see that only the voltages of $\sim 6$ nodes cross the limit by $0.005$ pu, the rest of the nodes generally maintain the same voltage level as the C-OPF results. Thus, we conclude that the proposed ENDiCO CLF
can maintain the voltages within the pre-specified limits even for an extreme DER penetration scenario.

\subsubsection{Solving Time}
In this section, we compare the required average solving time for C-OPF, previously developed ENDiCO controller, and the modified ENDiCO CLF controllers in Table \ref{time_tab}. The quantity represents the average time required to solve the OPF problem for each minute window. From the table we can see that C-OPF, ENDiCO with optimization solvers and ENDiCO with closed-form solutions take $11.4, 0.06$ and $2.6\times 10^{-5}$ s on an average, respectively to solve the VVC case. It is clear that using closed-form expressions with the ENDiCO controller greatly reduces the solving time.

\begin{small}
\begin{table}[h]
    \centering
    \caption{Solving Time Comparison}
    \vspace{-0.1cm}
        \label{time_tab}
    \setlength{\tabcolsep}{5pt}
    {\begin{tabular}{|c|c|c|c|}
    \hline
    \textbf{Cases} &C-OPF & ENDiCO & ENDiCO CLF\\
    \hline
    {\textbf{VVC}}& $11.4$ s & $0.06$ s& $2.6\times 10^{-5}$ s \\
    \hline
    {\textbf{VWC}}& $15.53$ s & $0.063$ s& $6.4\times 10^{-5}$ s \\
    \hline
    \end{tabular}}{}
   \vspace{-0.3cm}
\end{table}
\end{small}

\begin{small}
\begin{table}[h]
    \centering
    \caption{Controller Performance Summary}
    \vspace{-0.1cm}
        \label{Sum_tab1}
    \setlength{\tabcolsep}{5pt}
    {\begin{tabular}{|c|c|c|}
    \hline
    \textbf{Cases} &VVC &VWC\\
    \hline
    {\textbf{Max time steps to converge}}&2&4\\
    \hline
    {\textbf{Max voltage oscillation}}& 0.001 pu & 0.005 pu\\
    \hline
    {\textbf{Max tracking Error}}& 1.5\% & 1.8\%\\
    \hline
    {\textbf{Voltage limits Violation}}& 0 pu & 0.005 pu\\
    \hline
    \end{tabular}}{}
   \vspace{-0.3cm}
\end{table}
\end{small}

\subsection{Discussions}
The previously developed ENDiCO controller solves one of the major drawbacks of the state-of-the-art feedback-based online distributed controllers. It significantly reduces the number of time-steps required to track the optimal solutions, thus allowing a fast-tracking of optimal system operations with highly variables DERs. In this paper, the computational complexity of the controller is further reduced by incorporating the closed-form solutions with the ENDiCO controllers. This not only drastically minimizes the OPF solve time but also removes the necessity of placing an optimization solver at each distributed agent. Hence, the proposed approach significantly advances the state-of-the-art distributed feedback-based online voltage control applied to radial distribution systems. The computational performances are summarized in Table \ref{Sum_tab1}.

\section{Conclusions}
The proposed ENDiCO CLF controller leverages the radial topology of a power distribution system, significantly improving the computational time needed to obtain network-level optimal solutions.  The key improvement is the implementation of closed-form solutions to solve the node-level subproblems. This eliminates the need for embedded solvers at each agent, replacing them with a relatively simple algebraic computation. This new controller can handle rapid changes in system variables and accurately track optimal conditions even in highly stressed systems with minimal computation and communication requirements.  Solutions show minimal to no oscillations in nodal voltages.  In this paper, we considered a balanced network for algorithm design and simulations and we plan to expand the proposed method for an unbalanced three-phase distribution system.
\balance
\bibliographystyle{ieeetr}
\bibliography{cite}

\begin{thebibliography}{10}

\bibitem{momoh1999review1}
J.~A. Momoh, R.~Adapa, and M.~El-Hawary, ``A review of selected optimal power
  flow literature to 1993. i. nonlinear and quadratic programming approaches,''
  {\em IEEE transactions on power systems}, vol.~14, no.~1, pp.~96--104, 1999.

\bibitem{castillo2013survey}
A.~Castillo and R.~P. O’Neill, ``{Survey of approaches to solving the
  ACOPF},'' {\em US Federal Energy Regulatory Commission, Tech. Rep}, 2013.

\bibitem{erseghe2014distributed}
T.~Erseghe, ``Distributed optimal power flow using admm,'' {\em IEEE
  transactions on power systems}, vol.~29, no.~5, pp.~2370--2380, 2014.

\bibitem{dall2013distributed}
E.~Dall'Anese, H.~Zhu, and G.~B. Giannakis, ``Distributed optimal power flow
  for smart microgrids,'' {\em IEEE Transactions on Smart Grid}, vol.~4, no.~3,
  pp.~1464--1475, 2013.

\bibitem{millar2016smart}
B.~Millar and D.~Jiang, ``Smart grid optimization through asynchronous,
  distributed primal dual iterations,'' {\em IEEE Transactions on Smart Grid},
  vol.~8, no.~5, pp.~2324--2331, 2016.

\bibitem{magnusson2015distributed}
S.~Magn{\'u}sson, P.~C. Weeraddana, and C.~Fischione, ``A distributed approach
  for the optimal power-flow problem based on admm and sequential convex
  approximations,'' {\em IEEE Transactions on Control of Network Systems},
  vol.~2, no.~3, pp.~238--253, 2015.

\bibitem{bolognani2014distributed}
S.~Bolognani, R.~Carli, G.~Cavraro, and S.~Zampieri, ``Distributed reactive
  power feedback control for voltage regulation and loss minimization,'' {\em
  IEEE Transactions on Automatic Control}, vol.~60, no.~4, pp.~966--981, 2014.

\bibitem{cavraro2017local}
G.~Cavraro and R.~Carli, ``Local and distributed voltage control algorithms in
  distribution networks,'' {\em IEEE Transactions on Power Systems}, vol.~33,
  no.~2, pp.~1420--1430, 2017.

\bibitem{bernstein2019real}
A.~Bernstein and E.~Dall’Anese, ``Real-time feedback-based optimization of
  distribution grids: A unified approach,'' {\em IEEE Transactions on Control
  of Network Systems}, vol.~6, no.~3, pp.~1197--1209, 2019.

\bibitem{bastianello2020distributed}
N.~Bastianello, A.~Ajalloeian, and E.~Dall'Anese, ``Distributed and inexact
  proximal gradient method for online convex optimization,'' {\em arXiv
  preprint arXiv:2001.00870}, 2020.

\bibitem{qu2019optimal}
G.~Qu and N.~Li, ``Optimal distributed feedback voltage control under limited
  reactive power,'' {\em IEEE Transactions on Power Systems}, vol.~35, no.~1,
  pp.~315--331, 2019.

\bibitem{hu2019branch}
X.~Hu, Z.-W. Liu, G.~Wen, X.~Yu, and C.~Li, ``Branch-wise parallel successive
  algorithm for online voltage regulation in distribution networks,'' {\em IEEE
  Transactions on Smart Grid}, vol.~10, no.~6, pp.~6678--6689, 2019.

\bibitem{sadnan2020real}
R.~Sadnan and A.~Dubey, ``Real-time distributed control of smart inverters for
  network-level optimization,'' in {\em 2020 IEEE International Conference on
  Communications, Control, and Computing Technologies for Smart Grids
  (SmartGridComm)}, pp.~1--6, IEEE, 2020.

\bibitem{sadnan2021distributed}
R.~Sadnan and A.~Dubey, ``Distributed optimization using reduced network
  equivalents for radial power distribution systems,'' {\em IEEE Transactions
  on Power Systems}, 2021.

\bibitem{baran1989optimal2}
M.~Baran and F.~F. Wu, ``Optimal sizing of capacitors placed on a radial
  distribution system,'' {\em IEEE Transactions on power Delivery}, vol.~4,
  no.~1, pp.~735--743, 1989.

\bibitem{magnusson2020distributed}
S.~Magn{\'u}sson, G.~Qu, and N.~Li, ``Distributed optimal voltage control with
  asynchronous and delayed communication,'' {\em IEEE Transactions on Smart
  Grid}, 2020.

\end{thebibliography}

\end{document}